\documentclass[authoryear,12pt]{elsarticle}
\usepackage{cancel}
\usepackage{caption}
\usepackage{color}
\usepackage{lscape}
\usepackage{afterpage}
\usepackage{pstricks}
\usepackage{pst-plot}
\usepackage{longtable}
\usepackage{dcolumn}
\usepackage{pst-node}
\usepackage{amsmath}
\usepackage{amssymb}
\usepackage{amsthm}
\usepackage{multirow}
\usepackage{colortab}
\usepackage{color}
\usepackage{array}
\usepackage{slashbox}
\usepackage{colortbl}
\usepackage{subfigure}
\usepackage{textcomp}
\usepackage{pstricks, pst-node}
\usepackage{pst-all}
\usepackage{algorithm,algorithmic}
\usepackage{url}
\usepackage{soul}
\usepackage{hyperref}
\psset{arrows=->, labelsep=3pt, mnode=circle}

\usepackage{tabularx}

\usepackage{eurosym}
\usepackage{tikz}
\newcommand*\circled[1]{\tikz[baseline=(char.base)]{
\node[shape=circle,draw,inner sep=1pt] (char) {#1};}}

\newfont{\rams}{msbm10 scaled\magstep1}

\newenvironment{resumeT}{\begin{list}{}{\setlength{\rightmargin}{\leftmargin}}\item[]
{\centering {\bf \it~~~}
\par}\item[]\ignorespaces}{\unskip\end{list}}

\begin{document}
\title{Supporting the robust ordinal regression approach to multiple criteria decision aiding with a set of representative value functions}

\author[Eco]{\rm Sally Giuseppe Arcidiacono}
\ead{s.arcidiacono@unict.it}
\author[Eco]{\rm Salvatore Corrente}
\ead{salvatore.corrente@unict.it}
\author[Eco,por]{\rm Salvatore Greco}
\ead{salgreco@unict.it}

\address[Eco]{Department of Economics and Business, University of Catania, Corso Italia, 55, 95129  Catania, Italy}
\address[por]{University of Portsmouth, Portsmouth Business School, Centre of Operations Research and Logistics (CORL), Richmond Building, Portland Street, Portsmouth PO1 3DE, United Kingdom}

\date{}
\maketitle

\vspace{-1cm}

\begin{resumeT}

\textbf{Abstract:} \noindent In this paper we propose a new methodology to represent the results of the robust ordinal regression approach by means of a family of representative value functions for which, taken two alternatives $a$ and $b$, the following two conditions are satisfied: 1) if for all compatible value functions $a$ is evaluated not worse than $b$ and for at least one value function $a$ has a better evaluation, then the evaluation of $a$ is greater than the evaluation of $b$ for all representative value functions; 2) if there exists one compatible value function giving $a$ an evaluation greater than $b$ and another compatible value function giving $a$ an evaluation smaller than $b$, then there are also at least one representative function giving a better evaluation to $a$ and another representative value function giving $a$ an evaluation smaller than $b$. This family of representative value functions intends to provide the Decision Maker (DM) a more clear idea of the preferences obtained by the compatible value functions, with the aim to support the discussion in constructive approach of Multiple Criteria Decision Aiding.  

\end{resumeT}

 \pagenumbering{arabic}

\section{Introduction}

In this paper we extend and generalize the idea of representative value function of robust ordinal regression \cite{greco2011selection,kadzinski2012selection}. With respect to robust ordinal regression, in the following we consider the basic concepts and the notation in \cite{corrente2014robust}. The procedure to obtain an ``optimal'' family of representative value functions is based on the procedure permitting to find a compatible additive value  function $U$ such that
\begin{itemize}
	\item for all $a,b \in A$ such that $a \succ^N b$, that is $a \succsim^N b$ but not $b \succsim^N a$, we have $U(a) > U(b)$,
	\item for as many as possible $a,b \in A$ in a set $D \subseteq \{(a,b)\in A \times A: \mbox{ not } (a \succsim^N b) \mbox{ and not } (b \succsim^N a)\}$  we have $U(a) > U(b)$.
\end{itemize}

This function $U$ can be obtained solving the following MILP problem $P(D)$:
$$
min \sum_{(a,b)\in D}\gamma_{(a,b)}
$$
subject to
$$
\left.
\begin{array}{l}
\; U(a)=\sum_{i \in I}u_i(g_i(a))
\mbox{ for all $a\in A$}\\
\; u_i(g_i(a)) \ge u_i(g_i(b))
\mbox{ for all $a,b\in A$ such that $g_i(a) \ge g_i(b)$}\\
\; u_i(\alpha_i)=0,
\mbox{ for all $i\in I$}\\
\; \sum_{i \in I}u_i(\beta_i)=1\\
\;U(a)-U(b) \ge \varepsilon
\mbox{ if $a\succ b$}\\
\;U(a)-U(b) = 0
\mbox{ if $a \sim b$}\\
\;U(a)-U(b) \ge \varepsilon
\mbox{ if $a\succ^N b$}\\
\;U(a)-U(b) \ge \varepsilon -\gamma_{a,b}M
\mbox{ if $(a,b)\in D$}\\
\gamma_{(a,b)}\in\{0,1\}
\mbox{ for $(a,b)\in D$}\\
\end{array}\\
\right\}{U}_{Rep}^{D}
$$
\noindent where
\begin{itemize}
	\item  $\beta_j$ and $\alpha_j$ are the best and the worst considered values of criterion $g_j, j=1, \ldots, m$,
	\item  $\varepsilon$ is a small positive value, for example $\varepsilon=10^{-4}$,
	\item $M$ is a ``big number'', for example $M=10^4$.
\end{itemize}

To determine a ''sufficient set'' of representative value functions we can execute the following \textbf{Procedure 1}:
\begin{itemize}
	\item[Step 1.] Compute the necessary  preference relation $\succsim^N$;
	\item[Step 2.] Set $D=\{(a,b)\in A \times A: \mbox{ not } (a \succsim^N b) \mbox{ and not } (b \succsim^N a)\}$ and $r=1$;
	\item[Step 3.] Solve problem $P(D)$, denote by $U_r$ the value functions solution of the MILP problem, set $r^\prime=r$;
	\item[Step 4.] Set $D^\prime=D \setminus \{(a,b)\in A \times A: \gamma_{(a,b)}=0\}$;
	\item[Step 5.] If $D^\prime=\emptyset$, then $\mathcal{U}$=$\{U_1,\ldots,U_r\}$ is a ''sufficient set'' of representative value functions and the procedure terminates, otherwise set $D=D^\prime$, $r=r^\prime+1$ and go to Step 3.
\end{itemize}
Once with the previous procedure we have found a ''sufficient set'' of representative value functions, we can compute a minimal discriminant representative set of value functions by solving two MILP problems.

With the first MILP problem we compute the minimal number of value functions necessary to represent the relations $\succ^N$ and $\bowtie$ for which, for all $a,b \in A, a\bowtie b$ if not $(a \succsim^N b)$ and not $(b \succsim^N a)$. 
In particular, the following MILP problem $P1$ has to be solved
$$
max \sum_{s=1}^r\rho_s
$$
subject to
$$
\left.
\begin{array}{l}
\; U^s(a)=\sum_{i \in I}u^s_i(g_i(a))
\mbox{ for all $a\in A$}, s=1,\ldots,r\\
\; u^s_i(g_i(a)) \ge u^s_i(g_i(b))
\mbox{ for all $a,b\in A$ such that $g_i(a) \ge g_i(b), s=1,\ldots,r$}\\
\; u^s_i(\alpha_i)=0,
\mbox{ for all $i\in I, s=1,\ldots,r$}\\
\; \sum_{i \in I}u^s_i(\beta_i)=1, s=1,\ldots,r\\
\;U^s(a)-U^s(b) \ge \varepsilon
\mbox{ if $a\succ b$}\\
\;U^s(a)-U^s(b) = 0 
\mbox{ if $a\sim b$}\\
\;U^s(a)-U^s(b) \ge \varepsilon
\mbox{ if $a\succ^N b$}\\
\;U^s(a)-U^s(b) \ge \varepsilon-\gamma^s_{a,b}M-\rho^sM
\mbox{ if $a\bowtie b, s=1,\ldots,r$}\\
\; \gamma^s_{a,b} \le 1 -\rho^s, s=1,\ldots,r\\
\; \gamma^s_{(a,b)}\in\{0,1\}
\mbox{ for $a\bowtie b, s=1,\ldots,r$}\\
\; \sum_{s=1}^r \gamma^s_{a,b} \le r -\sum_{s=1}^r \rho^s -1 \mbox{ if $a\bowtie b$}\\
\; \rho^s \in \{0,1\}, s=1,\ldots,r\\
\end{array}\\
\right\}{U}_{Rep}^{Multi}
$$
\noindent where
\begin{itemize}
	\item  $\beta_j$ and $\alpha_j$, $\varepsilon$ and $M$ are as for the above problem $P(D)$,
	\item $r$ is the cardinal of the set $\mathcal{U}$ obtained with \textbf{Procedure 1}.
\end{itemize}

The minimal numbers of value functions $U$ necessary to represent the binary relations $\succ^N$ and $\bowtie$ is $t=r-z^*$ with $z^*=min \sum_{s=1}^r\rho_s$ obtained solving the MILP problem $P1$. Finally the most discriminant set of representative value functions $\mathcal{U}^{Discr}=\{U_1, \ldots, U_t\}$ is obtained solving the following MILP problem $P2$.    
$$
max \varepsilon
$$
subject to
$$
\left.
\begin{array}{l}
\; U^s(a)=\sum_{i \in I}u^s_i(g_i(a))
\mbox{ for all $a\in A$}, s=1,\ldots,t\\
\; u^s_i(g_i(a)) \ge u^s_i(g_i(b))
\mbox{ for all $a,b\in A$ such that $g_i(a) \ge g_i(b), s=1,\ldots,t$}\\
\; u^s_i(\alpha_i)=0,
\mbox{ for all $i\in I, s=1,\ldots,r$}\\
\; \sum_{i \in I}u^s_i(\beta_i)=1, s=1,\ldots,t\\
\;U^s(a)-U^s(b) \ge \varepsilon
\mbox{ if $a\succ b$}\\
\;U^s(a)-U^s(b) = 0 
\mbox{ if $a\sim b$}\\
\;U^s(a)-U^s(b) \ge \varepsilon
\mbox{ if $a\succ^N b$}\\
\;U^s(a)-U^s(b) \ge \varepsilon-\gamma^s_{a,b}M
\mbox{ if $a\bowtie b, s=1,\ldots,t$}\\
\; \gamma^s_{(a,b)}\in\{0,1\}
\mbox{ for $a\bowtie b, s=1,\ldots,t$}\\
\; \sum_{s=1}^r \gamma^s_{a,b} \le t -1 \mbox{ if $a\bowtie b$}\\
\end{array}\\
\right\}{U}_{Rep}^{Multi-Discr}
$$

\section{A didactic example}
We consider the problem, proposed by the Economist Intelligence Unit \cite{unit2010democracy}, of evaluating democracy of different countries on the basis of the following five criteria:
\begin{itemize}
	\item electoral process and pluralism (g1),
	\item the functioning of government (g2),
	\item political participation (g3),
	\item political culture (g4),
	\item civil liberties (g5).
\end{itemize}
Among the 165 countries in the ranking of \cite{unit2010democracy}, we considered the following ten countries: Indonesia, Japan, Kyrgyzstan, Malayasia, Mongolia, Papua New Guinea, Philippines, Singapore, Taiwan, Timor Leste. The performance table, with the evaluations of each country on the considered criteria, is given in Table \ref{tab:Perf_Table}.

\begin{table}[htb!]
  \centering
    \begin{tabular}{ccccccc}\hline
$a$ & {\tt Country} &$g_1$ & $g_2$ & $g_3$ & $g_4$ & $g_5$   \\ \hline
  $	a1	$&	Indonesia	(ID)	&	6.92	&	7.14	&	5.00	&	6.25	&	6.76	\\
$	a2	$&	Japan	(JP)	&	9.17	&	7.86	&	5.56	&	8.75	&	9.41	\\
$	a3	$&	Kyrgyzstan	(KG)	&	5.75	&	1.86	&	2.78	&	5.00	&	5.00	\\
$	a4	$&	Malayasia	(MY)	&	6.08	&	5.71	&	4.44	&	7.50	&	6.18	\\
$	a5	$&	Mongolia	(MN)	&	9.17	&	6.08	&	3.89	&	5.63	&	8.24	\\
$	a6	$&	Papua New Guinea 	(PG)	&	7.33	&	6.43	&	4.44	&	6.25	&	8.24	\\
$	a7	$&	Philippines	(PH)	&	9.17	&	5.36	&	5.00	&	3.75	&	9.12	\\
$	a8	$&	Singapore	(SG)	&	4.33	&	7.50	&	2.78	&	7.50	&	7.35	\\
$	a9	$&	Taiwan	(TW)	&	9.58	&	7.50	&	6.67	&	5.63	&	9.71	\\
$	a10	$&	Timor Leste	(TL)	&	7.00	&	5.57	&	5.00	&	6.25	&	8.24	\\  \hline
    \end{tabular}
  \caption{Performance table}\label{tab:Perf_Table}
\end{table}

We consider as Decision Maker (DM) an analyst providing the following preferences:
\begin{itemize}
	\item Malayasia is preffered to Mongolia, that is, $a4\succ a5$,
\item Singapore is preferred to Timor Leste, that is, $a8 \succ a10$,
\item Philippines is preferred to Papua New Guinea, that is, $a7 \succ a6$.
\end{itemize}
Applying the roubust ordinal regression \cite{greco2008ordinal}, we obtain the necessary preference relation represented in Table \ref{tab:NP_Table}.

\begin{table}[htb!]
  \centering
    \begin{tabular}{ccccccccccc}\hline
		&$	a1	$&$	a2	$&$	a3	$&$	a4	$&$	a5	$&$	a6	$&$	a7	$&$	a8	$&$	a9	$&$	a10	$\\ \hline
$	a1	$&	$\succsim^N$	&		&	$\succsim^N$	&		&		&		&		&		&		&		\\
$	a2	$&	$\succsim^N$	&	$\succsim^N$	&	$\succsim^N$	&	$\succsim^N$	&	$\succsim^N$	&	$\succsim^N$	&	$\succsim^N$	&	$\succsim^N$	&		&	$\succsim^N$	\\
$	a3	$&		&		&	$\succsim^N$	&		&		&		&		&		&		&		\\
$	a4	$&		&		&	$\succsim^N$	&	$\succsim^N$	&	$\succsim^N$	&		&		&		&		&		\\
$	a5	$&		&		&	$\succsim^N$	&		&	$\succsim^N$	&		&		&		&		&		\\
$	a6	$&		&		&	$\succsim^N$	&		&		&	$\succsim^N$	&		&		&		&		\\
$	a7	$&		&		&	$\succsim^N$	&		&		&	$\succsim^N$	&	$\succsim^N$	&		&		&		\\
$	a8	$&		&		&	$\succsim^N$	&		&		&		&		&	$\succsim^N$	&		&	$\succsim^N$	\\
$	a9	$&		&		&	$\succsim^N$	&		&	$\succsim^N$	&	$\succsim^N$	&	$\succsim^N$	&		&	$\succsim^N$	&		\\ 
$	a10	$&		&		&	$\succsim^N$	&		&		&		&		&		&		&	$\succsim^N$	\\ \hline
\end{tabular}
  \caption{Necessary preference $\succsim^N$}\label{tab:NP_Table}
\end{table}

On the basis of necessary preference relation $\succsim^N$, we can define the strict necessary preference relation $\succ^N$ and the relation $\bowtie$ shown in Tables \ref{tab:SNP_Table} and \ref{tab:INC_Table}, respectively.

\begin{table}[htb!]
  \centering
    \begin{tabular}{ccccccccccc}\hline
		&$	a1	$&$	a2	$&$	a3	$&$	a4	$&$	a5	$&$	a6	$&$	a7	$&$	a8	$&$	a9	$&$	a10	$\\ \hline
		$	a1	$&		&		&	$\succ^N$	&		&		&		&		&		&		&		\\
$	a2	$&	$\succ^N$	&		&	$\succ^N$	&	$\succ^N$	&	$\succ^N$	&	$\succ^N$	&	$\succ^N$	&	$\succ^N$	&		&	$\succ^N$	\\
$	a3	$&		&		&		&		&		&		&		&		&		&		\\
$	a4	$&		&		&	$\succ^N$	&		&	$\succ^N$	&		&		&		&		&		\\
$	a5	$&		&		&	$\succ^N$	&		&		&		&		&		&		&		\\
$	a6	$&		&		&	$\succ^N$	&		&		&		&		&		&		&		\\
$	a7	$&		&		&	$\succ^N$	&		&		&	$\succ^N$	&		&		&		&		\\
$	a8	$&		&		&	$\succ^N$	&		&		&		&		&	$\succ^N$	&		&	$\succ^N$	\\
$	a9	$&		&		&	$\succ^N$	&		&	$\succ^N$	&	$\succ^N$	&	$\succ^N$	&		&		&		\\
$	a10	$&		&		&	$\succ^N$	&		&		&		&		&		&		&		\\ \hline
\end{tabular}
  \caption{Strict necessary preference $\succ^N$}\label{tab:SNP_Table}
\end{table}

\begin{table}[htb!]
  \centering
    \begin{tabular}{ccccccccccc}\hline
		&$	a1	$&$	a2	$&$	a3	$&$	a4	$&$	a5	$&$	a6	$&$	a7	$&$	a8	$&$	a9	$&$	a10	$\\ \hline
		$	a1	$&		&		&		&	$\bowtie$	&	$\bowtie$	&	$\bowtie$	&	$\bowtie$	&	$\bowtie$	&	$\bowtie$	&	$\bowtie$	\\
$	a2	$&		&		&		&	&		&		&		&		&	$\bowtie$	&		\\
$	a3	$&		&		&		&		&		&		&		&		&		&		\\
$	a4	$&	$\bowtie$	&		&		&		&		&	$\bowtie$	&	$\bowtie$	&	$\bowtie$	&	$\bowtie$	&	$\bowtie$	\\
$	a5	$&	$\bowtie$	&		&		&		&		&	$\bowtie$	&	$\bowtie$	&	$\bowtie$	&		&	$\bowtie$	\\
$	a6	$&	$\bowtie$	&		&		&	$\bowtie$	&	$\bowtie$	&		&		&	$\bowtie$	&		&	$\bowtie$	\\
$	a7	$&	$\bowtie$	&		&		&	$\bowtie$	&	$\bowtie$	&		&		&	$\bowtie$	&		&	$\bowtie$	\\
$	a8	$&	$\bowtie$	&		&		&	$\bowtie$	&	$\bowtie$	&	$\bowtie$	&	$\bowtie$	&		&	$\bowtie$	&		\\
$	a9	$&	$\bowtie$	&	$\bowtie$	&		&	$\bowtie$	&		&		&		&	$\bowtie$	&		&	$\bowtie$	\\
$	a10	$&	$\bowtie$	&		&		&	$\bowtie$	&	$\bowtie$	&	$\bowtie$	&$\bowtie$&		&	$\bowtie$	&		\\\hline
\end{tabular}
  \caption{Incomparability $\bowtie$}\label{tab:INC_Table}
\end{table}

Using \textbf{Procedure 1}, the value functions 
\[
U^s(a)=\sum_{j=1}^5 u^s_j(g_j(a)), \;\; s=1,\ldots,4,
\]
presented in Tables \ref{tab:value_1}-\ref{tab:value_4} were found. They are sufficient to represent the binary relations $\succ^N$ and $\bowtie$.

\begin{table}[htb!]
  \centering
    \begin{tabular}{|cc|cc|cc|cc|cc|}\hline
		$g_1(a)$ & $u^1_1(g_1(a))$ & $g_2(a)$ &$u^1_2(g_2(a))$&$g_3(a)$ &$u^1_3(g_3(a))$&$g_4(a)$ &$u^1_4(g_4(a))$&$g_5(a)$ &$u^1_5(g_5(a))$\\ \hline
		4.33	&	0	&	1.86	&	0	&	2.78	&	0	&	3.75	&	0	&	5	&	0	\\
5.75	&	0	&	5.36	&	0	&	3.89	&	0	&	5	&	0	&	6.18	&	0	\\
6.08	&	0	&	5.57	&	0	&	4.44	&	0	&	5.63	&	0	&	6.76	&	0	\\
6.92	&	0	&	5.71	&	0	&	5	&	0.1667	&	6.25	&	0	&	7.35	&	0	\\
7	&	0.1666	&	6.08	&	0	&	5.56	&	0.1667	&	7.5	&	0.3331	&	8.24	&	0	\\
7.33	&	0.1666	&	6.43	&	0	&	6.67	&	0.1667	&	8.75	&	0.3331	&	9.12	&	0	\\
9.17	&	0.1668	&	7.14	&	0	&		&		&		&		&	9.41	&	0	\\
9.58	&	0.4998	&	7.5	&	0.0004	&		&		&		&		&	9.71	&	0	\\
	&		&	7.86	&	0.0004	&		&		&		&		&		&		\\ \hline
\end{tabular}
  \caption{Value function $U^1$}\label{tab:value_1}
\end{table}

\begin{table}[htb!]
  \centering
    \begin{tabular}{|cc|cc|cc|cc|cc|}\hline
		$g_1(a)$ & $u^2_1(g_1(a))$ & $g_2(a)$ &$u^2_2(g_2(a))$&$g_3(a)$ &$u^2_3(g_3(a))$&$g_4(a)$ &$u^2_4(g_4(a))$&$g_5(a)$ &$u^2_5(g_5(a))$\\ \hline
		4.33	&	0	&	1.86	&	0	&	2.78	&	0	&	3.75	&	0	&	5	&	0	\\
5.75	&	0	&	5.36	&	0	&	3.89	&	0	&	5	&	0	&	6.18	&	0	\\
6.08	&	0	&	5.57	&	0	&	4.44	&	0.166633333	&	5.63	&	0	&	6.76	&	0	\\
6.92	&	0	&	5.71	&	0.333366667	&	5	&	0.333266667	&	6.25	&	0	&	7.35	&	0.166433333	\\
7	&	0	&	6.08	&	0.333366667	&	5.56	&	0.333266667	&	7.5	&	0	&	8.24	&	0.166433333	\\
7.33	&	0	&	6.43	&	0.333366667	&	6.67	&	0.333266667	&	8.75	&	0	&	9.12	&	0.333266667	\\
9.17	&	0	&	7.14	&	0.333366667	&		&		&		&		&	9.41	&	0.333266667	\\
9.58	&	0	&	7.5	&	0.333366667	&		&		&		&		&	9.71	&	0.333266667	\\
	&		&	7.86	&	0.333366667	&		&		&		&		&		&		\\  \hline
\end{tabular}
  \caption{Value function $U^2$}\label{tab:value_2}
\end{table}

\begin{table}[htb!]
  \centering
    \begin{tabular}{|cc|cc|cc|cc|cc|}\hline
		$g_1(a)$ & $u^3_1(g_1(a))$ & $g_2(a)$ &$u^3_2(g_2(a))$&$g_3(a)$ &$u^3_3(g_3(a))$&$g_4(a)$ &$u^3_4(g_4(a))$&$g_5(a)$ &$u^3_5(g_5(a))$\\ \hline
		4.33	&	0	&	1.86	&	0	&	2.78	&	0	&	3.75	&	0	&	5	&	0	\\
5.75	&	0	&	5.36	&	0	&	3.89	&	0	&	5	&	0	&	6.18	&	0	\\
6.08	&	0	&	5.57	&	0.0001	&	4.44	&	0	&	5.63	&	0	&	6.76	&	0	\\
6.92	&	0	&	5.71	&	0.0001	&	5	&	0.374975	&	6.25	&	0.24985	&	7.35	&	0	\\
7	&	0	&	6.08	&	0.0001	&	5.56	&	0.374975	&	7.5	&	0.5	&	8.24	&	0	\\
7.33	&	0	&	6.43	&	0.125025	&	6.67	&	0.374975	&	8.75	&	0.5	&	9.12	&	0	\\
9.17	&	0	&	7.14	&	0.125025	&		&		&		&		&	9.41	&	0	\\
9.58	&	0	&	7.5	&	0.125025	&		&		&		&		&	9.71	&	0	\\
	&		&	7.86	&	0.125025	&		&		&		&		&		&		\\  \hline
\end{tabular}
  \caption{Value function $U^3$}\label{tab:value_3}
\end{table}

\begin{table}[htb!]
  \centering
    \begin{tabular}{|cc|cc|cc|cc|cc|}\hline
		$g_1(a)$ & $u^4_1(g_1(a))$ & $g_2(a)$ &$u^4_2(g_2(a))$&$g_3(a)$ &$u^4_3(g_3(a))$&$g_4(a)$ &$u^4_4(g_4(a))$&$g_5(a)$ &$u^4_5(g_5(a))$\\ \hline
		4.33	&	0	&	1.86	&	0	&	2.78	&	0	&	3.75	&	0	&	5	&	0	\\
5.75	&	0	&	5.36	&	0	&	3.89	&	0	&	5	&	0	&	6.18	&	0	\\
6.08	&	0	&	5.57	&	0.0001	&	4.44	&	0	&	5.63	&	0	&	6.76	&	0	\\
6.92	&	0	&	5.71	&	0.0001	&	5	&	0	&	6.25	&	0	&	7.35	&	0	\\
7	&	0	&	6.08	&	0.0001	&	5.56	&	0	&	7.5	&	0.0004	&	8.24	&	0	\\
7.33	&	0.0001	&	6.43	&	0.0001	&	6.67	&	0	&	8.75	&	0.0004	&	9.12	&	0	\\
9.17	&	0.0003	&	7.14	&	0.0001	&		&		&		&		&	9.41	&	0	\\
9.58	&	0.9995	&	7.5	&	0.0001	&		&		&		&		&	9.71	&	0	\\
	&		&	7.86	&	0.0001	&		&		&		&		&		&		\\ \hline
\end{tabular}
  \caption{Value function $U^4$}\label{tab:value_4}
\end{table}

Knowing that four value functions are sufficient to represent binary relations $\succ^N$ and $\bowtie$, one can solve MILP problem P1. So doing, it appears that three value functions are enough to represent binary relations $\succ^N$ and $\bowtie$. In particular, the three value functions $U^5, U^6$ and $U^7$ presented in Tables \ref{tab:value_5}-\ref{tab:value_7} were obtained from the solution of  MILP problem P1.  

\begin{table}[htb!]
  \centering
    \begin{tabular}{|cc|cc|cc|cc|cc|}\hline
		$g_1(a)$ & $u^5_1(g_1(a))$ & $g_2(a)$ &$u^5_2(g_2(a))$&$g_3(a)$ &$u^5_3(g_3(a))$&$g_4(a)$ &$u^5_4(g_4(a))$&$g_5(a)$ &$u^5_5(g_5(a))$\\ \hline
		4.33	&	0	&	1.86	&	0	&	2.78	&	0	&	3.75	&	0	&	5	&	0	\\
5.75	&	0	&	5.36	&	0	&	3.89	&	0.0005	&	5	&	0	&	6.18	&	0	\\
6.08	&	0	&	5.57	&	0	&	4.44	&	0.0005	&	5.63	&	0	&	6.76	&	0	\\
6.92	&	0	&	5.71	&	0.4995	&	5	&	0.0007	&	6.25	&	0	&	7.35	&	0	\\
7	&	0	&	6.08	&	0.4995	&	5.56	&	0.0007	&	7.5	&	0.0004	&	8.24	&	0	\\
7.33	&	0	&	6.43	&	0.4995	&	6.67	&	0.0007	&	8.75	&	0.0004	&	9.12	&	0.4991	\\
9.17	&	0.0003	&	7.14	&	0.4995	&		&		&		&		&	9.41	&	0.4991	\\
9.58	&	0.0003	&	7.5	&	0.4995	&		&		&		&		&	9.71	&	0.4991	\\
	&		&	7.86	&	0.4995	&		&		&		&		&		&		\\ \hline
\end{tabular}
  \caption{Value function $U^5$}\label{tab:value_5}
\end{table}

\begin{table}[htb!]
  \centering
    \begin{tabular}{|cc|cc|cc|cc|cc|}\hline
		$g_1(a)$ & $u^6_1(g_1(a))$ & $g_2(a)$ &$u^6_2(g_2(a))$&$g_3(a)$ &$u^6_3(g_3(a))$&$g_4(a)$ &$u^6_4(g_4(a))$&$g_5(a)$ &$u^6_5(g_5(a))$\\ \hline
		4.33	&	0	&	1.86	&	0	&	2.78	&	0	&	3.75	&	0	&	5	&	0	\\
5.75	&	0	&	5.36	&	0	&	3.89	&	0	&	5	&	0	&	6.18	&	0	\\
6.08	&	0	&	5.57	&	0	&	4.44	&	0.0005	&	5.63	&	0	&	6.76	&	0	\\
6.92	&	0	&	5.71	&	0	&	5	&	0.0006	&	6.25	&	0	&	7.35	&	0	\\
7	&	0	&	6.08	&	0	&	5.56	&	0.0006	&	7.5	&	0	&	8.24	&	0	\\
7.33	&	0.0004	&	6.43	&	0	&	6.67	&	0.0006	&	8.75	&	0	&	9.12	&	0	\\
9.17	&	0.0004	&	7.14	&	0.0002	&		&		&		&		&	9.41	&	0	\\
9.58	&	0.9987	&	7.5	&	0.0007	&		&		&		&		&	9.71	&	0	\\
	&		&	7.86	&	0.0007	&		&		&		&		&		&		\\ \hline
\end{tabular}
  \caption{Value function $U^6$}\label{tab:value_6}
\end{table}

\begin{table}[htb!]
  \centering
    \begin{tabular}{|cc|cc|cc|cc|cc|}\hline
		$g_1(a)$ & $u^7_1(g_1(a))$ & $g_2(a)$ &$u^7_2(g_2(a))$&$g_3(a)$ &$u^7_3(g_3(a))$&$g_4(a)$ &$u^7_4(g_4(a))$&$g_5(a)$ &$u^7_5(g_5(a))$\\ \hline
		4.33	&	0	&	1.86	&	0	&	2.78	&	0	&	3.75	&	0	&	5	&	0	\\
5.75	&	0	&	5.36	&	0	&	3.89	&	0	&	5	&	0	&	6.18	&	0	\\
6.08	&	0	&	5.57	&	0	&	4.44	&	0	&	5.63	&	0	&	6.76	&	0	\\
6.92	&	0	&	5.71	&	0	&	5	&	0.3333	&	6.25	&	0.3332	&	7.35	&	0.0001	\\
7	&	0	&	6.08	&	0	&	5.56	&	0.3333	&	7.5	&	0.6665	&	8.24	&	0.0001	\\
7.33	&	0	&	6.43	&	0	&	6.67	&	0.3333	&	8.75	&	0.6665	&	9.12	&	0.0001	\\
9.17	&	0	&	7.14	&	0	&		&		&		&		&	9.41	&	0.0001	\\
9.58	&	0	&	7.5	&	0	&		&		&		&		&	9.71	&	0.0001	\\
	&		&	7.86	&	0	&		&		&		&		&		&		\\   \hline
\end{tabular}
  \caption{Value function $U^7$}\label{tab:value_7}
\end{table}

Finally, knowing that three value functions are sufficient to represent binary relations $\succ^N$ and $\bowtie$, solving MILP problem P2, we can obtain the most discriminant set of representative value functions $\mathcal{U}^{Discr}=\{U^8, U^9, U^{10}\}$. The three value functions $U^8, U^9$ and $U^{10}$. corresponding to $\varepsilon^*=max \varepsilon=0,09091$, are presented in Tables \ref{tab:value_8}-\ref{tab:value_10}, while the pair $(a,b) \in A \times A$ represented by value functions $U^8, U^9$ and $U^{10}$ are shown in Tables \ref{tab:Rep_U_8}-\ref{tab:Rep_U_10}. The evaluations assigned by the marginal value functions and the overall evaluations assigned by value representative value functions $U^8, U^9$ and $U^{10}$ to alternatives (Countries) $a_1,\ldots,a_{10}$ are shown in Tables \ref{tab:Ev_U_8}-\ref{tab:Ev_U_10}. Observe that the pairs $(a,b) \in A \times A$ for which $a \succ^N b$ are represented by all the three value functions $U^8, U^9$ and $U^{10}$, while the pair $(a,b) \in A \times A$ for which $a \bowtie b$ are represented by at least one of the three value functions $U^8, U^9$ and $U^{10}$. In this perspective, the DM could be interested to see a compatible value function for which $U(a)>U(b)$ for some compatible value function $U$ and for this  $U^8, U^9$ and $U^{10}$ provide always at least one example of such a value function. For example, if the DM is interested to see a compatible value function for which $U(a4) > U(a8)$, the value function $U^{10}$ can be shown, so that, observing that the only marginal value functions assigning a different evaluation to $a4$ and $a8$ is $u^{10}_3(g_3(\cdot))$, one can conclude that an explanation for assigning $a4$ a better evaluation than $a8$ is in the performances assigned to the two alternatives by criterion $g_3$ political participation for which $a4$ (Malaysya) has a score of 4.44 and $a8$ (Singapore) a score of 2.78.

\begin{table}[htb!]
  \centering
    \begin{tabular}{|cc|cc|cc|cc|cc|}\hline
		$g_1(a)$ & $u^8_1(g_1(a))$ & $g_2(a)$ &$u^8_2(g_2(a))$&$g_3(a)$ &$u^8_3(g_3(a))$&$g_4(a)$ &$u^8_4(g_4(a))$&$g_5(a)$ &$u^8_5(g_5(a))$\\ \hline
		4.33	&	0	&	1.86	&	0	&	2.78	&	0	&	3.75	&	0	&	5	&	0	\\
5.75	&	0	&	5.36	&	0	&	3.89	&	0	&	5	&	0	&	6.18	&	0	\\
6.08	&	0	&	5.57	&	0.090909091	&	4.44	&	0	&	5.63	&	0	&	6.76	&	0	\\
6.92	&	0	&	5.71	&	0.090909091	&	5	&	0	&	6.25	&	0	&	7.35	&	0	\\
7	&	0.136363636	&	6.08	&	0.090909091	&	5.56	&	0	&	7.5	&	0.5	&	8.24	&	0	\\
7.33	&	0.136363636	&	6.43	&	0.090909091	&	6.67	&	0	&	8.75	&	0.5	&	9.12	&	0	\\
9.17	&	0.318181818	&	7.14	&	0.136363636	&		&		&		&		&	9.41	&	0	\\
9.58	&	0.318181818	&	7.5	&	0.181818182	&		&		&		&		&	9.71	&	0	\\
	&		&	7.86	&	0.181818182	&		&		&		&		&		&		\\   \hline
\end{tabular}
  \caption{Value function $U^8$}\label{tab:value_8}
\end{table}

\begin{table}[htb!]
  \centering
    \begin{tabular}{|cc|cc|cc|cc|cc|}\hline
		$g_1(a)$ & $u^9_1(g_1(a))$ & $g_2(a)$ &$u^9_2(g_2(a))$&$g_3(a)$ &$u^9_3(g_3(a))$&$g_4(a)$ &$u^9_4(g_4(a))$&$g_5(a)$ &$u^9_5(g_5(a))$\\ \hline
		4.33	&	0	&	1.86	&	0	&	2.78	&	0	&	3.75	&	0	&	5	&	0	\\
5.75	&	0	&	5.36	&	0	&	3.89	&	0	&	5	&	0	&	6.18	&	0	\\
6.08	&	0	&	5.57	&	0	&	4.44	&	0	&	5.63	&	0	&	6.76	&	0	\\
6.92	&	0	&	5.71	&	0	&	5	&	0.295454545	&	6.25	&	0.204545455	&	7.35	&	0.113636364	\\
7	&	0	&	6.08	&	0	&	5.56	&	0.295454545	&	7.5	&	0.5	&	8.24	&	0.113636364	\\
7.33	&	0	&	6.43	&	0	&	6.67	&	0.295454545	&	8.75	&	0.5	&	9.12	&	0.113636364	\\
9.17	&	0	&	7.14	&	0.090909091	&		&		&		&		&	9.41	&	0.113636364	\\
9.58	&	0	&	7.5	&	0.090909091	&		&		&		&		&	9.71	&	0.113636364	\\
	&		&	7.86	&	0.090909091	&		&		&		&		&		&		\\ \hline
\end{tabular}
  \caption{Value function $U^9$}\label{tab:value_9}
\end{table}

\begin{table}[htb!]
  \centering
    \begin{tabular}{|cc|cc|cc|cc|cc|}\hline
		$g_1(a)$ & $u^{10}_1(g_1(a))$ & $g_2(a)$ &$u^{10}_2(g_2(a))$&$g_3(a)$ &$u^{10}_3(g_3(a))$&$g_4(a)$ &$u^{10}_4(g_4(a))$&$g_5(a)$ &$u^{10}_5(g_5(a))$\\ \hline
		4.33	&	0	&	1.86	&	0	&	2.78	&	0	&	3.75	&	0	&	5	&	0	\\
5.75	&	0	&	5.36	&	0	&	3.89	&	0	&	5	&	0	&	6.18	&	0	\\
6.08	&	0	&	5.57	&	0	&	4.44	&	0.181818182	&	5.63	&	0	&	6.76	&	0	\\
6.92	&	0	&	5.71	&	0.272727273	&	5	&	0.181818182	&	6.25	&	0	&	7.35	&	0	\\
7	&	0	&	6.08	&	0.272727273	&	5.56	&	0.181818182	&	7.5	&	0	&	8.24	&	0	\\
7.33	&	0.090909091	&	6.43	&	0.272727273	&	6.67	&	0.272727273	&	8.75	&	0	&	9.12	&	0.363636364	\\
9.17	&	0.090909091	&	7.14	&	0.272727273	&		&		&		&		&	9.41	&	0.363636364	\\
9.58	&	0.090909091	&	7.5	&	0.272727273	&		&		&		&		&	9.71	&	0.363636364	\\
	&		&	7.86	&	0.272727273	&		&		&		&		&		&		\\ \hline
\end{tabular}
  \caption{Value function $U^{10}$}\label{tab:value_10}
\end{table}

\begin{table}[htb!]
  \centering
    \begin{tabular}{ccccccccccc}\hline
		&$	a1	$&$	a2	$&$	a3	$&$	a4	$&$	a5	$&$	a6	$&$	a7	$&$	a8	$&$	a9	$&$	a10	$\\ \hline
$	a1	$ &	 	&	 	&	\circled{$\succ^N$}	&	$\bowtie$	&	$\bowtie$	&	$\bowtie$	&	$\bowtie$	&	$\bowtie$	&	$\bowtie$	&	$\bowtie$	\\
$	a2	$ &	\circled{$\succ^N$}	&	 	&	\circled{$\succ^N$}	&	\circled{$\succ^N$}	&	\circled{$\succ^N$}	&	\circled{$\succ^N$}	&	\circled{$\succ^N$}	&	\circled{$\succ^N$}	&	\circled{$\bowtie$}	&	\circled{$\succ^N$}	\\
$	a3	$ &	 	&	 	&	 	&	 	&	 	&	 	&	 	&	 	&	 	&	 	\\
$	a4	$ &	\circled{$\bowtie$}	&	 	&	\circled{$\succ^N$}	&	 	&	\circled{$\succ^N$}	&	\circled{$\bowtie$}	&	\circled{$\bowtie$}	&	$\bowtie$	&	\circled{$\bowtie$}	&	\circled{$\bowtie$}	\\
$	a5	$ &	\circled{$\bowtie$}	&	 	&	\circled{$\succ^N$}	&	 	&	 	&	\circled{$\bowtie$}	&	\circled{$\bowtie$}	&	$\bowtie$	&	 	&	\circled{$\bowtie$}	\\
$	a6	$ &	\circled{$\bowtie$}	&	 	&	\circled{$\succ^N$}	&	$\bowtie$	&	$\bowtie$	&	 	&	 	&	$\bowtie$	&	 	&	$\bowtie$	\\
$	a7	$ &	\circled{$\bowtie$}	&	 	&	\circled{$\succ^N$}	&	$\bowtie$	&	$\bowtie$	&	\circled{$\succ^N$}	&	 	&	$\bowtie$	&	 	&	\circled{$\bowtie$}	\\
$	a8	$ &	\circled{$\bowtie$}	&	 	&	\circled{$\succ^N$}	&	\circled{$\bowtie$}	&	\circled{$\bowtie$}	&	\circled{$\bowtie$}	&	\circled{$\bowtie$}	&	 	&	\circled{$\bowtie$}	&	\circled{$\succ^N$}	\\
$	a9	$ &	\circled{$\bowtie$}	&	$\bowtie$	&	\circled{$\succ^N$}	&	$\bowtie$	&	\circled{$\succ^N$}	&	\circled{$\succ^N$}	&	\circled{$\succ^N$}	&	$\bowtie$	&	 	&	\circled{$\bowtie$}	\\
$	a10	$ &	\circled{$\bowtie$}	&	 	&	\circled{$\succ^N$}	&	$\bowtie$	&	$\bowtie$	&	$\bowtie$	&	$\bowtie$	&	 	&	$\bowtie$	&	 		\\ \hline
\end{tabular}
  \caption{Pairs $(a_h,a_k)$ of the relations $\succsim^N$ and $\bowtie$ represented by value function $U^8$}\label{tab:Rep_U_8}
\end{table}

\begin{table}[htb!]
  \centering
    \begin{tabular}{ccccccccccc}\hline
		&$	a1	$&$	a2	$&$	a3	$&$	a4	$&$	a5	$&$	a6	$&$	a7	$&$	a8	$&$	a9	$&$	a10	$\\ \hline
$	a1	$ &	 	&	 	&	\circled{$\succ^N$}	&	\circled{$\bowtie$}	&	\circled{$\bowtie$}	&	\circled{$\bowtie$}	&	\circled{$\bowtie$}	&	$\bowtie$	&	\circled{$\bowtie$}	&	$\bowtie$	\\
$	a2	$ &	\circled{$\succ^N$}	&	 	&	\circled{$\succ^N$}	&	\circled{$\succ^N$}	&	\circled{$\succ^N$}	&	\circled{$\succ^N$}	&	\circled{$\succ^N$}	&	\circled{$\succ^N$}	&	\circled{$\bowtie$}	&	\circled{$\succ^N$}	\\
$	a3	$ &	 	&	 	&	 	&	 	&	 	&	 	&	 	&	 	&	 	&	 	\\
$	a4	$ &	$\bowtie$	&	 	&	\circled{$\succ^N$}	&	 	&	\circled{$\succ^N$}	&	\circled{$\bowtie$}	&	$\bowtie$	&	$\bowtie$	&	$\bowtie$	&	$\bowtie$	\\
$	a5	$ &	$\bowtie$	&	 	&	\circled{$\succ^N$}	&	 	&	 	&	$\bowtie$	&	$\bowtie$	&	$\bowtie$	&	 	&	$\bowtie$	\\
$	a6	$ &	$\bowtie$	&	 	&	\circled{$\succ^N$}	&	$\bowtie$	&	\circled{$\bowtie$}	&	 	&	 	&	$\bowtie$	&	 	&	$\bowtie$	\\
$	a7	$ &	$\bowtie$	&	 	&	\circled{$\succ^N$}	&	$\bowtie$	&	\circled{$\bowtie$}	&	\circled{$\succ^N$}	&	 	&	$\bowtie$	&	 	&	$\bowtie$	\\
$	a8	$ &	\circled{$\bowtie$}	&	 	&	\circled{$\succ^N$}	&	\circled{$\bowtie$}	&	\circled{$\bowtie$}	&	\circled{$\bowtie$}	&	\circled{$\bowtie$}	&	 	&	\circled{$\bowtie$}	&	\circled{$\succ^N$}	\\
$	a9	$ &	$\bowtie$	&	$\bowtie$	&	\circled{$\succ^N$}	&	$\bowtie$	&	\circled{$\succ^N$}	&	\circled{$\succ^N$}	&	\circled{$\succ^N$}	&	$\bowtie$	&	 	&	$\bowtie$	\\
$	a10	$ &	$\bowtie$	&	 	&	\circled{$\succ^N$}	&	\circled{$\bowtie$}	&	\circled{$\bowtie$}	&	\circled{$\bowtie$}	&	\circled{$\bowtie$}	&	 	&	\circled{$\bowtie$}	&	 		\\ \hline
\end{tabular}
  \caption{Pairs $(a_h,a_k)$ of the relations $\succsim^N$ and $\bowtie$ represented by value function $U^9$}\label{tab:Rep_U_9}
\end{table}

\begin{table}[htb!]
  \centering
    \begin{tabular}{ccccccccccc}\hline
		&$	a1	$&$	a2	$&$	a3	$&$	a4	$&$	a5	$&$	a6	$&$	a7	$&$	a8	$&$	a9	$&$	a10	$\\ \hline
$	a1	$ &	 	&	 	&	\circled{$\succ^N$}	&	$\bowtie$	&	\circled{$\bowtie$}	&	$\bowtie$	&	$\bowtie$	&	\circled{$\bowtie$}	&	$\bowtie$	&	\circled{$\bowtie$}	\\
$	a2	$ &	\circled{$\succ^N$}	&	 	&	\circled{$\succ^N$}	&	\circled{$\succ^N$}	&	\circled{$\succ^N$}	&	\circled{$\succ^N$}	&	\circled{$\succ^N$}	&	\circled{$\succ^N$}	&	$\bowtie$	&	\circled{$\succ^N$}	\\
$	a3	$ &	 	&	 	&	 	&	 	&	 	&	 	&	 	&	 	&	 	&	 	\\
$	a4	$ &	$\bowtie$	&	 	&	\circled{$\succ^N$}	&	 	&	\circled{$\succ^N$}	&	$\bowtie$	&	$\bowtie$	&	\circled{$\bowtie$}	&	$\bowtie$	&	\circled{$\bowtie$}	\\
$	a5	$ &	$\bowtie$	&	 	&	\circled{$\succ^N$}	&	 	&	 	&	$\bowtie$	&	$\bowtie$	&	\circled{$\bowtie$}	&	 	&	\circled{$\bowtie$}	\\
$	a6	$ &	\circled{$\bowtie$}	&	 	&	\circled{$\succ^N$}	&	\circled{$\bowtie$}	&	\circled{$\bowtie$}	&	 	&	 	&	\circled{$\bowtie$}	&	 	&	\circled{$\bowtie$}	\\
$	a7	$ &	\circled{$\bowtie$}	&	 	&	\circled{$\succ^N$}	&	\circled{$\bowtie$}	&	\circled{$\bowtie$}	&	\circled{$\succ^N$}	&	 	&	\circled{$\bowtie$}	&	 	&	\circled{$\bowtie$}	\\
$	a8	$ &	$\bowtie$	&	 	&	\circled{$\succ^N$}	&	$\bowtie$	&	$\bowtie$	&	$\bowtie$	&	$\bowtie$	&	 	&	$\bowtie$	&	\circled{$\succ^N$}	\\
$	a9	$ &	\circled{$\bowtie$}	&	\circled{$\bowtie$}	&	\circled{$\succ^N$}	&	\circled{$\bowtie$}	&	\circled{$\succ^N$}	&	\circled{$\succ^N$}	&	\circled{$\succ^N$}	&	$\bowtie$	&	 	&	\circled{$\bowtie$}	\\
$	a10	$ &	$\bowtie$	&	 	&	\circled{$\succ^N$}	&	$\bowtie$	&	$\bowtie$	&	$\bowtie$	&	$\bowtie$	&	 	&	$\bowtie$	&	 		\\ \hline
\end{tabular}
  \caption{Pairs $(a_h,a_k)$ of the relations $\succsim^N$ and $\bowtie$ represented by value function $U^{10}$}\label{tab:Rep_U_10}
\end{table}

\begin{table}[htb!]
  \centering
    \begin{tabular}{|c|c|c|c|c|c|c|c|}\hline
	& Country	&$	u^8_1(g_1(a))$	&$	u^8_2(g_2(a))	$&$	u^8_3(g_3(a))	$&$	u^8_4(g_4(a))	$&$	u^8_5(g_5(a))	$&$	U^8(a) $\\ \hline
$	a1	$&	Indonesia	&	0	&	0.136363636	&	0	&	0	&	0	&	0.136363636	\\
$	a2	$&	Japan	&	0.318181818	&	0.181818182	&	0	&	0.5	&	0	&	1	\\
$	a3	$&	Kyrgyzstan	&	0	&	0	&	0	&	0	&	0	&	0	\\
$	a4	$&	Malayasia	&	0	&	0.090909091	&	0	&	0.5	&	0	&	0.590909091	\\
$	a5	$&	Mongolia	&	0.318181818	&	0.090909091	&	0	&	0	&	0	&	0.409090909	\\
$	a6	$&	Papua New G. 	&	0.136363636	&	0.090909091	&	0	&	0	&	0	&	0.227272727	\\
$	a7	$&	Philippines	&	0.318181818	&	0	&	0	&	0	&	0	&	0.318181818	\\
$	a8	$&	Singapore	&	0	&	0.181818182	&	0	&	0.5	&	0	&	0.681818182	\\
$	a9	$&	Taiwan	&	0.318181818	&	0.181818182	&	0	&	0	&	0	&	0.5	\\
$	a10	$&	Timor Leste	&	0.136363636	&	0.090909091	&	0	&	0	&	0	&	0.227272727	\\ \hline
\end{tabular}
  \caption{Evaluations assigned by the marginal value functions and overall evaluations with respect to value function $U^8$}\label{tab:Ev_U_8}
\end{table}

\begin{table}[htb!]
  \centering
    \begin{tabular}{|c|c|c|c|c|c|c|c|}\hline
	& Country	&$	u^8_1(g_1(a))$	&$	u^8_2(g_2(a))	$&$	u^8_3(g_3(a))	$&$	u^8_4(g_4(a))	$&$	u^8_5(g_5(a))	$&$	U^8(a) $\\ \hline
$	a1	$&	Indonesia	&	0	&	0.090909091	&	0.295454545	&	0.204545455	&	0	&	0.590909091	\\
$	a2	$&	Japan	&	0	&	0.090909091	&	0.295454545	&	0.5	&	0.113636364	&	1	\\
$	a3	$&	Kyrgyzstan	&	0	&	0	&	0	&	0	&	0	&	0	\\
$	a4	$&	Malayasia	&	0	&	0	&	0	&	0.5	&	0	&	0.5	\\
$	a5	$&	Mongolia	&	0	&	0	&	0	&	0	&	0.113636364	&	0.113636364	\\
$	a6	$&	Papua New G. 	&	0	&	0	&	0	&	0.204545455	&	0.113636364	&	0.318181818	\\
$	a7	$&	Philippines	&	0	&	0	&	0.295454545	&	0	&	0.113636364	&	0.409090909	\\
$	a8	$&	Singapore	&	0	&	0.090909091	&	0	&	0.5	&	0.113636364	&	0.704545455	\\
$	a9	$&	Taiwan	&	0	&	0.090909091	&	0.295454545	&	0	&	0.113636364	&	0.5	\\
$	a10	$&	Timor Leste	&	0	&	0	&	0.295454545	&	0.204545455	&	0.113636364	&	0.613636364	\\
 \hline
\end{tabular}
  \caption{Evaluations assigned by the marginal value functions and overall evaluations with respect to value function $U^{9}$}\label{tab:Ev_U_9}
\end{table}

\begin{table}[htb!]
  \centering
    \begin{tabular}{|c|c|c|c|c|c|c|c|}\hline
	& Country	&$	u^8_1(g_1(a))$	&$	u^8_2(g_2(a))	$&$	u^8_3(g_3(a))	$&$	u^8_4(g_4(a))	$&$	u^8_5(g_5(a))	$&$	U^8(a) $\\ \hline
$	a1	$&	Indonesia	&	0	&	0.272727273	&	0.181818182	&	0	&	0	&	0.454545455	\\
$	a2	$&	Japan	&	0.090909091	&	0.272727273	&	0.181818182	&	0	&	0.363636364	&	0.909090909	\\
$	a3	$&	Kyrgyzstan	&	0	&	0	&	0	&	0	&	0	&	0	\\
$	a4	$&	Malayasia	&	0	&	0.272727273	&	0.181818182	&	0	&	0	&	0.454545455	\\
$	a5	$&	Mongolia	&	0.090909091	&	0.272727273	&	0	&	0	&	0	&	0.363636364	\\
$	a6	$&	Papua New G. 	&	0.090909091	&	0.272727273	&	0.181818182	&	0	&	0	&	0.545454545	\\
$	a7	$&	Philippines	&	0.090909091	&	0	&	0.181818182	&	0	&	0.363636364	&	0.636363636	\\
$	a8	$&	Singapore	&	0	&	0.272727273	&	0	&	0	&	0	&	0.272727273	\\
$	a9	$&	Taiwan	&	0.090909091	&	0.272727273	&	0.272727273	&	0	&	0.363636364	&	1	\\
$	a10	$&	Timor Leste	&	0	&	0	&	0.181818182	&	0	&	0	&	0.181818182	\\
\hline
\end{tabular}
  \caption{Evaluations assigned by the marginal value functions and overall evaluations with respect to value function $U^{10}$}\label{tab:Ev_U_10}
\end{table}

\section{Conclusions}

We presented a methodology to represent the results of the robust ordinal regression approach through a set of representative value functions in multiple criteria decision aiding procedures. We plan to extend the idea of a set of representative models also to robust ordinal regression approach applied to sorting problems and to outranking models. We plan also to apply our methodology in real world decision problems to test its effective support.

\bibliographystyle{plainnat}
\bibliography{Full_bibliography}


\end{document}